\documentclass[12pt]{article}
\usepackage{theorem}
\usepackage{amssymb}
\usepackage{amsmath}

\newtheorem{thm}{Theorem}[section]
\newtheorem{prop}[thm]{Proposition}
\newtheorem{cor}[thm]{Corollary}
\newtheorem{lm}[thm]{Lemma}

{\theorembodyfont{\rm}
\newtheorem{exe}[thm]{Exercise}
\newtheorem{defn}[thm]{Definition}
\newtheorem{exam}[thm]{Example}
\newtheorem{rem}[thm]{Remark}

}

\def\ra{\rightarrow}
\def\lg{{\langle}}
\def\rg{{\rangle}}
\def\ind{{\rm ind}}

\def\C{{\mathbb C}}

\def\Ga{{\Gamma}}
\def\H{{\rm H}}
\def\N{{\mathbb N}}
\def\M{{\bf M}}
\def\P{{\mathbb P}}
\def\Q{{\mathbb Q}}
\def\R{{\mathbb R}}

\def\Z{{\mathbb Z}}

\def\Mod{{\bf M}}
\def\L{{\Lambda_{\rm eff}}}
\def\Pic{{\rm Pic}}
\def\rk{{\rm rk}}
\def\Aut{{\rm Aut}}
\def\o{{\omega}}
\def\s{{\sigma}}
\def\dim{{\rm dim}}

\def\SL{{\rm SL}}
\def\a{{\alpha}}

\title{Density of rational points on \\
elliptic K3 surfaces}

\author{F. A.  Bogomolov\\
\small  Courant Institute of Mathematical Sciences, N.Y.U. \\
\small 251 Mercer str. \\
\small New York, (NY) 10012, U.S.A.\\
\small e-mail: bogomolo@cims.nyu.edu\\
\small  and\\
Yu. Tschinkel\\
\small Dept. of Mathematics, U.I.C.\\
\small 851 South Morgan str.\\
\small Chicago, (IL) 60607-7045,  U.S.A.  \\
\small e-mail: yuri@math.uic.edu
}

\begin{document}

\date{}

\maketitle

\thispagestyle{empty}

\pagebreak

\section{Introduction}

Let $X$ be a smooth projective algebraic variety defined
over a number field $K$. We will say that rational points
on $X$ are potentially dense if there exists a finite extension
$K'/K$ such that the set $X(K')$ of  
$K'$-rational points is Zariski dense.
What are possible strategies to propagate rational points
on an algebraic variety? We thought of two:
using the group of automorphisms $\Aut(X)$ and
using additional geometric structures - like elliptic 
fibrations. The class of K3 surfaces 
is an ideal test case for both methods.

\noindent
One of our main results is:

\begin{thm}\label{thm:main-thm}
Let $X$ be a K3 surface defined over a number field $K$.
Assume that $X$ has a structure of an elliptic fibration
or an infinite group of automorphisms. 
Then rational points on $X$ are potentially dense. 
\end{thm}
Here is a more detailed list of what we learned:
We don't know if rational points are potentially dense on 
a general K3 surface with Picard group $\Pic(X_{\C})=\Z$.
In particular, we don't know if rational points are dense
on a double cover of 
$\P^2 $ ramified in a general curve of degree 6.
However, we can prove potential density for a {\em divisor} in 
the space of all such K3 surfaces, corresponding to the case
when, for example, the ramification curve is singular 
(cf. \cite{bogomolov-tschi-98-2}). 
The overall picture is similar. 
In any moduli family of algebraic K3 surfaces
we can find some union of algebraic subsets, 
including a divisor, such 
that rational points are potentially dense on the K3 surfaces 
corresponding to the points of this subset. More precisely, 

\begin{thm}
Let $X$ be a K3 surface, defined over a number field $K$.
Assume that  $\rk\,\Pic(X_{\C})=2$ and that 
$X$ does not contain a $(-2)$-curve.
Then rational points on $X$ are potentially dense. 
\end{thm}

\begin{rem}
If $\rk\,\Pic(X_{\C})=2$ and if $X$ does not 
contain a $(-2)$-curve then 
either it has an elliptic fibration or it has
an infinite automorphism group (but not both!).
For example, a quartic surface in 
$\P^3$ containing a smooth curve 
of genus 2 and degree 6 doesn't admit any elliptic fibrations, but
the group $\Aut(X_{\C})$ is infinite 
(cf. \cite{P-Sha} p. 583, \cite{Severi}). 
\end{rem}

\begin{thm} Let $X$ be a 
K3 surface over $K$ with ${\rm rk}\,\Pic(X_{\C})\ge 3$. 
Then rational points on $X$ are potentially dense, with a possible
exception of 8 isomorphy classes of lattices $\Pic(X_{\C})$.
\end{thm}

\begin{rem}\label{rem:aut}
If ${\rm rk}\,\Pic(X_{\C})= 3$ then there are only 6 types of lattices where
we can't prove potential density. There are only 2 types when 
${\rm rk}\Pic(X_{\C})= 4$. 
Potential density holds for all K3 with 
${\rm rk}\Pic(X_{\C})\ge 5$. 
All K3 surfaces with $\rk \, \Pic(X_{\C})=20$ have
infinite groups of automorphisms. 
We use Nikulin's classification 
of lattices of algebraic K3 surfaces 
(cf. \cite{nikulin-86}, \cite{nikulin-81}). 
\end{rem}

First we consider the problem of density for 
general elliptic fibrations ${\cal E}\ra \P^1$. 
Suppose that ${\cal E}$ has a zero section 
(i.e. ${\cal E}$ is Jacobian) and that there exists a section
of infinite order in the Mordell-Weil group of ${\cal E}$. 
Then a specialization argument shows
that rational points are dense in ${\cal E}_b$ 
for a Zariski dense set of
fibers $b\in \P^1$ (cf. \cite{silverman}). 
It turns out that even in absence of
global sections one can sometimes
arrive at the same conclusion. 

\begin{defn}\label{dfn:nt}
Let ${\cal E}\ra B$ be an elliptic fibration and 
${\cal M}\subset {\cal E}$ 
an irreducible multisection (defined over $\C$)
with the following property:  
for a general point $b\in B(\C)$ there exist two distinct points 
$p_b,p_b'\in ({\cal M}\cap {\cal E}_b)(\C)$ such that
$p_b-p_b'$ is non-torsion in the 
Jacobian ${\cal J}({\cal E}_b)(\C)$ of ${\cal E}_b$. 
We will call such a multisection an $nt$-multisection 
(non-torsion).
\end{defn}

For example, if ${\cal M}$ is ramified in a 
smooth fiber of ${\cal E}$ then it is an $nt$-multisection
(cf. \ref{prop:salient-non-torsion}). 
We will say that ${\cal M}$ is {\em torsion of order $m$} if 
for all $b\in B$ and 
all $p_b,p_b'\in {\cal M}\cap {\cal E}_b$ the zero-cycle
$p_b-p_b'$ is torsion of order $m$ in ${\cal J}({\cal E}_b)$.
An easy lemma (but not a tautology!) says that if ${\cal M}$ 
is not torsion of order $m$ for any $m\in \N$ then ${\cal M}$ is
an $nt$-multisection (cf. \ref{lm:nt}). 
(There are analogous notions for abelian schemes 
and torsors under abelian schemes.) 

\begin{prop}\label{prop:non-torsion}
Assume that ${\cal E}\ra \P^1$ has an
$nt$-multisection which is a 
rational or elliptic curve.
Then rational points on 
${\cal E}$ are potentially dense.
\end{prop}

We want to study situations when rational or
elliptic multisections occur and to analyze 
constrains which they impose on the elliptic
fibration (possible monodromy, structure
of singular fibers etc). 
We shall call fibrations with finitely many (resp. none)
rational or elliptic multisections 
{\em hyperbolic} (resp. {\em strongly hyperbolic}). 
Unfortunately,  we don't know examples of 
hyperbolic elliptic fibrations (without multiple fibers).
The aim of Section 2 is to prove the existence of a least
{\em one} rational multisection on algebraic elliptic K3 surfaces. 
From this we will deduce in Section 3 
the following theorem:

\begin{thm}\label{thm:nt-inf}
Let $X$ be an algebraic K3 surface 
with $\rk \, \Pic(X_{\C})\le 19$ admitting 
a structure of an elliptic fibration. 
Then this fibration has infinitely 
many rational $nt$-multisections. 
\end{thm}

The proof goes roughly as follows:
We find elliptic K3 surfaces ${\cal E}'\ra \P^1$
admitting a dominant map ${\cal E}'\ra X$ such that
the genus of every irreducible $m$-torsion 
multisection ${\cal M}'\subset {\cal E}'$ is $\ge 2$. 
On the other hand, the deformation theory argument 
in Section 2
implies that ${\cal E}'$ contains a {\em rational} multisection 
which must be an $nt$-multisection. Its image in $X$ is
a rational $nt$-multisection. 

\

{\small 
{\bf Acknowledgements.} 
Part of this work was done while 
both authors were visiting the Max-Planck-Institute
in Bonn. We are grateful to the MPI for the hospitality. 
The first author was partially  supported by the NSF. The second
author was partially supported by the NSA.
We would like to thank Joe Harris and Barry Mazur for their
ideas, suggestions and encouragement.
}

\section{K3 surfaces}

In this section we prove that every elliptic fibration
on an algebraic K3 surface has at least
one rational multisection. 

\subsection{Generalities}

There are several approaches to the theory of K3 surfaces.
Algebraically, a K3 surface $S$ (defined over some field of
characterstic zero) is a smooth projective
surface with trivial canonical class ${\cal K}_S=0$  
and $\H^1(S,{\cal O}_S)=0$. They are parametrized by an infinite
countable set of 19-dimensional algebraic spaces. 
The main invariant is the Picard group $\Pic(S)$ which is
isomorphic to a torsion free primitive lattice 
of finite rank ($\le 20$) equipped with a hyperbolic even integral
bilinear form. 

\

Another approach is via K\"ahler geometry. 
A K3 surface $S$ is a  compact simply connected  K\"ahler surface
equipped with a non-degenerate nowhere vanishing  holomorphic $(2,0)$-form
$\omega_S$. To obtain a natural parametrization we have to consider
{\em marked} K3 surfaces, which are pairs $(S,\s)$ consisting
of a K3 surface $S$ and an isometry of lattices
$$
\s : \H^2(S,\Z)\stackrel{\sim}{\ra} {\cal L}\simeq
3\cdot{\cal H} \oplus 2\cdot(- E_8),
$$ 
where ${\cal H} $ is the standard lattice with form $xy$ and
$E_8$ is an 8-dimensional even unimodular positive definite lattice.
We will denote by $\lg , \rg$ the intersection form 
on $\H^2(S,\Z)$. 
Marked K3 surfaces are parametrized by the conformal class
$\s_{\C}(\H^{2,0}(S,\Z))$ of their
non-degenerate holomorphic forms - the {\em period}. The latter
lies in the quadric given by  $\lg \o_S,\o_S\rg =0$ 
(inside $\P^{21}= \P(\H^2(S,\Z)_{\C})$). The period (still denoted by)
$\o_S$ satisfies the inequality $\lg \o_s,\overline{\o}_S\rg >0$.
Therefore, marked K\"ahler K3 surfaces are parametrized by
points of a complex homogeneous domain 
$\Omega = SO_{(3,19)}(\R)/SO_{(2,18)}(\R)$. (with the standard equivariant
complex structure). Unmarked K3 surfaces correspond to orbits
of the group $SO_{(3,19)}(\Z)$ on this space.

\

We will identify cycles and forms on $S$ with their (co)homology
classes. We will call a homology class $h$ {\em primitive}, if
$h\neq m Z$ for some $m>1$ and some effective cycle $Z$. 
We denote by $\L(S)$ the monoid of all classes in $\Pic(S)$ 
represented by effective divisors. (This differs slightly from the standard
definition of the effective cone as a cone in $\Pic(S)_{\R}$. 
In particular,  the smallest closed  {\em cone } in $\Pic(S)_{\R}$ 
containing $\L(S)$ could be finitely generated with  
$\L(S)$  being infinitely generated.)

\

We want to describe, in this setting, the subset of algebraic
and elliptic K3 surfaces. A K\"ahler K3 surface  $S$ is algebraic if
there is a primitive element $x\in \H^2(S,\Z)$ such that $\lg x,x\rg >0$
and $\lg \o_S,x\rg =0$. 
Conversely, every primitive $x\in \H^2(S,\Z)$  
determines a hyperplane $\{ \lg \o_S,x\rg =0\}$.
The intersection of this hyperplane with $\Omega$ will be
denoted by $\Omega(x)$. For a generic point of $\Omega(x)$ 
with $\lg x,x\rg >0$ one of the classes $\pm x$ defines a 
polarization of the corresponding 
marked K3 surface.

\

Every element $h$ which is a generator of $\L(S)$ 
with $\lg h,h\rg = -2$ is represented by a smooth rational curve.
Similarly, every generator 
of $\L(S)$ with $\lg h,h\rg = 0$ is represented
by a smooth elliptic curve (which defines an elliptic fibration 
without multiple fibers $S\ra \P^1$). In particular, this class is
also represented by a (singular) rational curve, contained in 
the singular fibers of the fibration. 
Therefore,  (marked) elliptic K3 
surfaces constitute a set of hyperplanes
$\Omega(h)$ with $\lg h,h\rg =0$ (and primitive $h$). 
For a generic member of $\Omega(h)$
the element $h$ defines the class of a fiber of the corresponding 
elliptic fibration.

\subsection{Deformation theory}

In this section  we work over $\C$. 
An {\em immersion} of a smooth curve $f: C\ra X$  
into a smooth variety $X$ is a regular map of degree 1 onto its
image such that the differential ${\rm d}f$ is non-zero everywhere. 
An {\em embedding} is an immersion with smooth image.

\begin{rem}
If $f : C\ra S$ is an immersion 
of a smooth curve into a smooth surface
then there exists a local neighborhood $U$ of $C$ (abstractly)
to which the map $f$ extends as a local isomorphism $f : U \ra S$.
The normal bundle ${\cal N}_C(U)$ of $C$ in $U$ is
defined by restriction of the canonical bundle ${\cal K}_S$ to
$f(C)$. In particular, if $S$ is a K3 surface then the normal bundle
${\cal N}_C(U)={\cal K}_C$.
\end{rem}

\begin{prop}\label{prop:deform-curve}
Let $C_0$ be a smooth rational curve, $S_0$  a K3 surface 
and $f_0 : C_0\ra S_0 $ an immersion.
Let  ${\cal S}\ra {\cal T}$ be a smooth scheme over a
complex ball ${\cal T}$ of dimension 20 with fibers
smooth K3 surfaces $S_t$ (local deformations of $S_0$). 
Consider the smooth subfamily 
${\cal S}'_{{\cal T}'}={\cal S}'\ra {\cal T}'$ corresponding
to deformations such that the class of 
$[f_0(C_0)]\in \H^2(S_0,\Z)\simeq\H^2(S_t,\Z)$ 
remains algebraic for all $t\in {\cal T}'$ (dimension of ${\cal T}'$ 
equals 19).  
Then for all $t\in {\cal T}'$ (close enough to $t_0$) 
there exists a smooth family of smooth curves 
${\cal C}_{{\cal T}'}= {\cal C}\ra {\cal T}' $ and a holomorphic map 
$f' : {\cal C}_{{\cal T}'}\ra  {\cal S}_{{\cal T}'}$  such that 
$f'|_{t_0}=f_0$ 
\end{prop}

{\em Proof.}
Construct a complex 2-dimensional neighborhood  $U_0$ of 
$C_0$ with the property that 
$f_0$ extends to a holomorphic map $g_0 : U_0\ra S_0$
such that $g_0$ is a local isomorphism. This is
possible since ${\rm d}f_0 \neq 0$. 
There is a non-degenerate $(2,0)$-form on $U_0$ 
induced from $S_0$. The curve $C_0$ is smooth in $U_0$ 
and its normal bundle in 
$U_0$ is isomorphic to ${\cal O}_{C_0}(-2)$.
It is well known that in this situation
there exists a local neighborhood of $C_0$ which is isomorphic to 
a small neighborhood of the zero section in the bundle ${\cal O}_{C_0}(-2)$.

The deformation of the complex structure on $S_0$ induces (by means
of $g_0$) a deformation of the complex structure on $U_0$. We 
obtain a smooth family $g : {\cal U}_{\cal T} \ra {\cal S}_{\cal T}$ 
(with $g|_{t_0}=g_0$) of deformations of
complex structures on $U_0$. 
The base of the space of versal deformations
for $U_0$ is a 1-dimensional disc. 
In the neighborhood
of $t_0\in {\cal T}$ every deformation of $U_0$ is induced
from the versal deformation space by a holomorphic map. 
As a preimage of zero we obtain a local divisor ${\cal D}_0\subset 
{\cal T}$. It follows that $C_0\times {\cal D}_0$ is contained
in the restriction of the family ${\cal U}_{\cal T}$ to ${\cal D}_0$.

On the other hand, outside the divisor ${\cal T}'\subset {\cal T}$
the class $[C_0]\in \H^2(U_t,\Z)=\Z$ is not algebraic. This is equivalent
to the property that the integral of the holomorphic form 
$\omega_t$ over the class $[f_0(C_0)]$ is not zero. 
Then the integral of the induced form 
$g^*(\omega_t)$ over $[C_0]$ is not equal to zero as well 
(where $\omega_t$ is the non-degenerate 
holomorphic form on $S_t$ induced by deformation). 
Therefore, the class $[C_0]\in \H^2(U_t,\Z)$ cannot be
realized by a holomorphic curve if $t\not\in {\cal T}'$.
Since we have obtained a 
realization of this class over ${\cal D}_0$ we can conclude that 
the local divisor  ${\cal D}_0$ is contained in  ${\cal T}'$.
Since ${\cal T}'$ is irreducible (it is a smooth disc), both 
divisors coincide.  Therefore, the map $f'$ is obtained by
restriction of  $g$ to ${\cal C}_{{\cal D}_0} = {\cal C}_{{\cal T}'}$.  

\begin{rem}
This proof imitates the approach of S. Bloch 
who introduced the notion of semi-regularity 
for {\em embedded} varieties (\cite{bloch}).
Here we use a similar technique for immersed varieties.
This deformation technique was extended to the case of
general maps by Z. Ran (cf. \cite{ran} and \cite{ran-2}).
\end{rem}

\subsection{Effective divisors}

\begin{thm}{\rm (Bogomolov-Mumford)}
Every class in $\L(S)$ can be
represented by a sum of (classes of) rational curves.
\end{thm}  

{\em Proof.}
The monoid of effective divisors $\L(S)$ of a K3 surface $S$ 
is generated by classes of $(-2)$-curves (represented
by smooth rational curves), classes $x$ with $\lg x,x\rg=0$
(represented by smooth elliptic curves, cuspidal elliptic curves
or nodal elliptic curves) and by 
primitive classes $x$ with $\lg x,x\rg>0$. 
Any smooth elliptic curve defines an elliptic fibration.
This fibration always has singular fibers (Euler characteristic)
and they consist of rational curves. 
It remains to show the following

\begin{prop} 
Let $S$ be a K3 surface. Every primitive effective class in $\Pic(S)$
with $\lg x,x\rg >0$ 
can be represented by a sum of (classes of) rational curves
(with multiplicities).
\end{prop}

The rest of this section is devoted to a proof of this
fact. An alternative proof is contained in \cite{mori-mukai-83-2}.

\

We first show that every primitive class is 
uniquely determined by its square.
Next we give a direct construction of a K3 surface
containing a rational curve which represents a primitive class,
with a given square. Finally we apply a deformation argument. 
 
\begin{exe} \label{exe}
Let ${\cal L}$ be an indefinite unimodular lattice
containing $ 3\cdot {\cal H}$, where ${\cal H}$
is the standard form given by $xy$. 
The orbit of any primitive element  
under the group ${\rm SO}({\cal L})$ is 
uniquely determined by the square
of this element. 
\end{exe}

{\em Proof.} First we show it for elements with $x$ with
$\lg x,x \rg=0$. Indeed, since $x$ is primitive, there exists a $y$ with
$\lg x, y\rg =1$. Then $x, y$ generate a sublattice ${\cal H}\subset {\cal L}$.
Since any sublattice ${\cal H}$ is a direct summand, we have 
the result. Similarly, if $z$ is any element such that $\lg z, x\rg = 1$
for some $x\in {\cal L}$ with $\lg x,x\rg =0$ then $z$ is equivalent
to the element with coordinates $(\lg z,z\rg,1)$ in the sublattice ${\cal H}$
(with standard coordinates). This concludes the exercise
(see also \cite{Isk-Shaf}, p. 224).

\begin{cor} Let $S$ be a K3 surface. 
Every primitive class in $\L(S)$ is 
uniquely determined by its self-intersection. 
\end{cor}

{\em Proof.} Identify $\Pic(S)$ with a sub-lattice in
${\cal L}= 3\cdot  {\cal H} \oplus 2\cdot (- E_8)$.

\begin{prop}
For any even $n\in \N$ there exists a pair $f : C\hookrightarrow S$
consisting of a smooth rational curve $C$ immersed
in a K3 surface $S$ and having self-intersection
equal to $n$. 
\end{prop}

{\em Proof.}
Let $R$ be a curve of genus 2 and ${\cal J}(R)$ its Jacobian. 
Let $\Z/\ell\Z \subset {\cal J}(R)$ 
be a cyclic subgroup of odd order $\ell$. 
Consider the map $\pi : R\ra {\cal J}(R)/(\Z/\ell\Z)$. 

\begin{lm}\label{lm:genericity}
For a generic $R$ the curve $\pi(R)$ contains exactly 6 points of order 2
of the quotient abelian variety  ${\cal J}(R)/(\Z/\ell\Z)$. 
These points are non-singular points of $\pi(R)$.
\end{lm}

{\em Proof.} 
It suffices to show that the only torsion points of ${\cal J}(R)$ contained
in $R$ (for a generic $R$) are the standard 6 points of order 2.
(Indeed, a point $\pi(Q)$, where $Q\in R$ is a point 
of order 2 in ${\cal J}(R)$ is a singular point of 
$\pi(R)$ if and only if there exists a point $P\neq Q$ in $R$ such that
$\ell \cdot P = Q$ in ${\cal J}(R)$. 
Thus $P$ has to be a torsion point of order
$2\ell$.)  

Consider the universal family ${\cal C}\ra \Mod(2,2)$ 
of smooth curves of genus $2$ with 2 level structure.
This family is imbedded as a subvariety into the
universal family of principally
polarized abelian varieties ${\cal J}\ra \Mod_J(2,2)$ (Jacobians) 
of dimension $2$ and level $2$. 
The family ${\cal C}\ra \Mod(2,2)$ has 6 natural sections
(points of order 2 in the Jacobian).
The family ${\cal J}\ra \Mod_J(2,2)$  has $16$ natural sections 
and $6$ of them are contained in ${\cal C}$.

The monodromy 
of the family ${\cal C}\ra \Mod(2,2)$ is 
a congruence subgroup 
of the group ${\rm Sp}_4(\Z)$ (which we
denote by $\Gamma_{\cal C}$). 
The torsion multisections of ${\cal J}\ra \Mod_J(2,2)$
split into a countable union of irreducible varieties $T_m$ corresponding 
to the orbits of monodromy $\Gamma_{\cal C}$ on $(\Q/\Z)^4$.  

Thus if a generic element $R_t$ contains a torsion 
point of order $\ell$ then it also contains its $\Gamma_{\cal C}$ orbit. 
Remark that for odd $\ell$ this orbit  consists of all 
primitive elements of order $\ell$ in 
the torsion group of the fiber. 
If $\ell = 2n$ then the corresponding orbit
contains all primitive torsion points $x,y$ of order $\ell$ 
with $nx = ny$.
Thus the intersection cycle $ R_t + a \cap   R_t$ 
(where $a$ is an element
of order $\ell $ or $n = \ell/2$ in the 
even case) consists of primitive points 
$x$ such that $x + a$ is also primitive.
The degree of this cycle is 
$\geq \phi(\ell) \ell^2$  (where $\phi(\ell)$ is the Euler function.) 
(In fact, for any primitive $a$ the degree 
is greater than the number of points
$x$ which are primitive modulo the subgroup generated by $a$ - hence 
the number of primitive points in 
$(\Z/\Z \ell)^3$ estimates the corresponding number
from below.) On the other hand $ \lg R_t + a , R_t\rg  = 2$.
Hence we obtain a contradiction if $m > 2$.

It shows that the only torsion points which can lie on a 
generic curve of genus $2$
are the points of order $2$.
Since any point of order $2$ which lies on $R_t$ has to be invariant
under the standard involution there are exactly six points of this kind 
on any $R$.

\begin{lm} 
The self-intersection $ \lg \pi(R),\pi(R) \rg=2\ell$.
\end{lm}

{\em Proof.} Indeed, the preimage $\pi^{-1}(\pi(R))$ 
consists of translations of $R$
by $\Z/\ell\Z$. 
Since $\lg R,R \rg =2$ we have 
$ \lg \pi(R),\pi(R) \rg =\frac{1}{\ell}\cdot 2\ell^2$.

\begin{lm}
For every even $n>0$ there exists a K3 surface $S_n$ containing
a rational curve which represents a primitive
class $c_n$ such that  $\lg c_n, c_n\rg=n$.
\end{lm}

{\em Proof.}
After dividing  ${\cal J}(R)/(\Z/\ell\Z)$ by 
$\Z/2\Z$ we obtain a rational curve $\pi(R)/(\Z/2\Z)$ 
on the singular Kummer surface ${\cal J}(R)/{\cal D}_{2\ell}$
(where ${\cal D}_{2\ell}$ is the dihedral group). 
After blowing up ${\cal J}(R)/{\cal D}_{2\ell}$ at the images
of the 16 points of order 2 on ${\cal J}(R)/(\Z/\ell\Z)$
we obtain an immersed rational curve with square
$\ell-3$. This curve represents a primitive  class
because its intersection with
each of the 6 blown up $(-2)$-curves equals to one.

\begin{lm} Let $S$ be any K3 surface with an
effective primitive class $x$ with square equal to $n$.
Then there exists a 1-dimensional smooth family
of K3 surfaces $f : {\cal S}\ra {\cal T}_1$ such 
that $x$ is an effective class in 
$\Pic(S_t)$ for all $t\in {\cal T}$ and
such that $S_{t_0}=S$, $S_{t_1}=S_n$  (for some $t_0,t_1\in {\cal T}$)
and the class $x\in \Pic(S_{t_1})$ is represented by 
an immersed rational curve. 
\end{lm}

{\em Proof.}
Consider a subvariety in the moduli space of marked K3 surfaces
where a given class $x$ is algebraic. 
It is given by a hyperplane section 
with the equation $\lg \o,x\rg =0$
in the intersection of the
open domain $\lg \o, \overline{\o}\rg >0$ 
with the quadric $\lg \o,\o\rg =0$. 
This is a connected smooth domain, which 
we denote by $\Omega(x)$. This domain is invariant under the action 
of a subgroup of $SO_{(2,19)}(\R)$. The arithmetic subgroup 
$\Gamma(x)\subset SO_{(3,19)}(\Z)$ stabilizing $x$ 
acts discretely on $\Omega(x)$ (since it stabilizes a
3-dimensional subspace generated by $x,\o,\overline{\o}$ which 
has a positive definite intersection form) and the quotient 
is a possibly singular 
algebraic variety with at most 
quotient singularities. It is a (coarse) moduli space 
of K3 surfaces with a fixed class $x$. 
There exists
a subgroup $\Gamma(x)'$ of finite index in $\Gamma(x)$  which acts
freely on $\Omega(x)$. The quotient $\tilde{\Omega}(x)$ is the fine
moduli space of K3 surfaces with a fixed class $x$ such that
for a generic point of $\tilde{\Omega}(x)$ the corresponding 
K3 surface carries a polarization 
with class $x$. There is a  point in $\Omega(x)$
which corresponds to a Kummer surface $S_n$ with a class
$c_n=x$ represented by an immersed rational curve. 
(Indeed, the classes $c_n,x$ lie in the same orbit under the action
of $SO_{(3,19)}(\Z)$.) 
We have a smooth algebraic curve
$\tilde{{\cal T}}(x) \subset\tilde{\Omega}(x)$
which connects the projections of points corresponding to $S_n$ and $S$.
Observe that we can choose the curve $\tilde{{\cal T}}(x)$ 
such that it contains only a 
finite number of points $\tilde{t}$ where $x$ is not
a polarization of the corresponding K3 surface $S_{\tilde{t}}$.
The family of  effective cycles  $C_{\tilde{t}}(x)$
(represented by sums  of rational curves) which represent the class $x$ in 
the group $\Pic(S_{\tilde{t}})$ is an algebraic ruled surface 
which projects surjectively onto the generic point of 
$\tilde{{\cal T}}(x)$ (this follows from the surjectivity in 
the neighborhood of $S_n$). 
Hence, there is a smooth relative compactification of this ruled 
surface with a proper (fiberwise) map to the corresponding
family of K3 surfaces. The  class of the image of any fiber 
$C_{\tilde{t}}(x)$ is $x$. 

\begin{rem}\label{rem:prim}
Let $x$ be a primitive class which 
is one of the generators of $\L(S)$. Then 
it is represented by an irreducible rational curve.
\end{rem}

\begin{rem}
There are similar results about immersions
of stable curves (not necessarily rational curves) and 
substantially more general theorems on the existence of
curves and families of curves. For example, 
Mori and Mukai proved that a generic K3 surface 
can be covered by a family of elliptic curves (cf. \cite{mori-mukai-83}).
Yau, Zaslow and Beauville found a formula for 
the number of (singular) rational curves in a given class on generic
K3 surfaces (\cite{yau-zaslow}, \cite{beauville-97}). 
Xi Chen constructs such 
curves deforming them from combinations of rational curves
on degenerations of K3 surfaces (cf. \cite{chen}, \cite{chen-2}).
(However, their results don't imply the existence of infinitely 
many rational multisections on elliptic K3 surfaces.) 
Let us also mention the work of C. Voisin on Lagrangian immersions
of algebraic varieties into hyperk\"ahler varieties. 
We decided to include
the initial argument of the first author since it is direct,  
transparent and sufficient for our purposes. 

\end{rem}

\begin{prop} 
The set $\M(S,h)$ of elliptic K3 surfaces ${\cal E}\ra \P^1$ 
with a fixed Jacobian ${\cal J}({\cal E})=S$ is
given by $\M(S,h)=\{ \omega_t=\omega_S+th\}_{t\in \C} \subset
\Omega(h)$,
where $t$ is a complex parameter and $h$ is a representative
of the class of the elliptic fiber $S_b$ ($b\in \P^1$). 
\end{prop}  

{\em Proof.} 
Let $h$ be any $(1,1)$-form induced from 
the base $\P^1$. Then the form $\o_h:=\o_S+th$ defines a complex
structure on $S$. Indeed, $\o_h$ is non-zero everywhere,
its square is identically zero, it is a closed form  and
it is non-degenerate on the real sub-bundle of the tangent 
bundle. If its class is homologous to zero then the
variation is trivial. Otherwise, we obtain a line $\M(S,h)$
in the space  $\Omega(h)$. 

Assume now that ${\cal E}'\ra \P^1$ is an elliptic K3 surface 
with the same given Jacobian $S$. 
Then there is a smooth (fiberwise) isomorphism 
$\iota : {\cal E}\ra {\cal E}'$ which is 
holomorphic along the fibers.  
The holomorphic forms $\o_{\cal E}$ and $\omega_{{\cal E}'}$ 
correspond to the sections $s,s'$ of $\H^0(\P^1,{\cal O})$.
Therefore, the difference  
$\omega_{\cal E} - \iota^*(\omega_{{\cal E}'})$ is a closed form
which has a non-trivial kernel on the tangent sub-bundle
to elliptic fibers. Therefore, this difference is a form of
rank at most 2 induced from the base of the elliptic fibration.

\section{Elliptic fibrations}

\subsection{Generalities}

In this section we continue to work over $\C$. 
We have to use parallel theories of elliptic fibrations
in the analytic and in the algebraic categories. 
All algebraic constructions carry over to the analytic
category. As in the case of K3 surfaces there are
some differences which we explain along the way.

\begin{defn}
Let ${\cal E}$ be a smooth projective algebraic surface.
An elliptic fibration is a morphism
$\varphi :  {\cal E}\ra B$ onto a smooth projective irreducible curve $B$ 
with connected fibers and with generic fiber a smooth curve of genus  $1$. 
A Jacobian elliptic fibration is an elliptic fibration with 
a section $e : B\ra {\cal E}$. 
\end{defn}

To every elliptic fibration $\varphi : {\cal E}\ra B$ one can 
associate a Jacobian elliptic fibration 
$\varphi_{\cal J}: {\cal J}={\cal J}({\cal E})\ra B$ 
(cf. \cite{barth-peters-vdv-84}). 
Over the generic point  ${\cal J}_b$ is given by classes of divisors
of degree zero in the fiber ${\cal E}_b$. 
The zero section $e_{\cal J}$ corresponds to the trivial class. 
A Jacobian elliptic fibration ${\cal J}$ can be viewed simultaneously as
a group scheme over $B$ (defining a sheaf over $B$) and as a surface
(the total space). In order to distinguish, we will sometimes use
the notation ${\cal J}$ and $S({\cal J})$, respectively. 
Most of the time we will work with $B=\P^1$.

\

We will only consider elliptic fibrations without multiple fibers.
They are locally isomorphic to the associated Jacobian elliptic
fibration ${\cal J}={\cal J}({\cal E})$ (for every point in the base $b\in B$
there exists a neighborhood $U_b$ such that the fibration
${\cal E}$ restricted to $U_b\subset B$ is Jacobian). 
The fibration ${\cal E}$ is a principal homogeneous space (torsor) under
${\cal J}$ and the set of all (isomorphism classes of)
${\cal E}$ with fixed Jacobian is identified with 
$\H^1_{\rm et}(B,{\cal J})$
(where ${\cal J}$ is considered as a sheaf of sections in the 
Jacobian elliptic fibration). 
In the analytic category we have a similar description of
elliptic fibrations ${\cal E}$ with given Jacobian
${\cal J}$ (where ${\cal J}$ is  always algebraic).
The group of isomorphism classes of ${\cal E}$ with a
given Jacobian ${\cal J}$ is identified with $\H^1_{\rm an}(B,{\cal J})$.

\

In the presence of singular fibers
we have 
$$
\H^1(B,{\cal J})= \H^2(S({\cal J}),{\cal O})/
{\rm Image}(\H^2(S({\cal J}),\Z)).
$$
The subgroup of algebraic elliptic
fibrations $\H^1_{\rm et}(B,{\cal J})$ 
coincides with the torsion subgroup in this quotient
(\cite{friedman-morgan-94}, Section 1.5). 
It can  also be described  as the union of
the images of  $\H^1_{\rm an}(B, {\cal J}_m)$, noting  
the exact sequence
$$
\H^0_{\rm an}(B, {\cal J})\ra \H^1_{\rm an}(B, {\cal J}_m)\ra 
\H^1_{\rm an}(B,{\cal J})
$$
where ${\cal J}_m$ is the sheaf of elements of order $m$ in ${\cal J}$
(the elements of order $m$ lie in the image of 
$\H^1_{\rm an}(B, {\cal J}_m)$).

\

\subsection{Multisections}

\begin{defn} Let 
$\varphi :  {\cal E}\ra B$ be an 
elliptic fibration (analytic or
algebraic).
We say  that a subvariety (analytic or algebraic)
${\cal M}\subset {\cal E}$ is a {\em multisection} of
degree $d_{\cal E}({\cal M})$  if 
${\cal M} $ is irreducible and if the degree $d_{\cal E}({\cal M})$ 
of the projection
$\varphi :  {\cal M}\ra B$ is non-zero. The definition of degree
extends to formal linear combinations of multisections.
\end{defn}

\begin{rem}
If an analytic fibration ${\cal E}\ra B$ has an analytic multisection 
then both the fibration and the multisection are algebraic.
\end{rem}

\noindent
There is a natural map 
$$
{\rm Rest}\,:\, \Pic({\cal E})\ra \Pic({\cal E}_b)/\Pic^0({\cal E}_b)=\Z.
$$

\begin{defn} The {\em degree}  $d_{\cal E}$ of
an algebraic elliptic fibration  ${\cal E}\ra B$ 
is the index of the image of $\Pic({\cal E})$ under the map
${\rm Rest}$. 
\end{defn}
Clearly, the degree of any multisection ${\cal M}$
of ${\cal E}$  is divisible by $d_{\cal E}$.

\begin{lm}
There exists a multisection ${\cal M}\subset {\cal E}$
with $d_{\cal E}({\cal M})=d_{\cal E}$. 
\end{lm}

{\em Proof.}
Let ${\cal D}$ be a divisor in ${\cal E}$ 
representing the class having intersection $d_{\cal E}$ with 
the class of the generic fiber of ${\cal E}$. 
Then there is an effective divisor in the class of 
${\cal D}'={\cal D}+ n\cdot{\cal E}_b$
for some $n\ge 0$. Indeed, consider 
$\lg {\cal D}',{\cal D}'\rg = 
\lg {\cal D},{\cal D}\rg +2n d_{\cal E}({\cal D})$.
By Riemann-Roch, the Euler characteristic 
is 
$$
\frac{1}{2}\lg {\cal D}',{\cal D}'-{\cal K}_{\cal E}\rg 
+c_1^2+\frac{c_2}{12}= \frac{1}{2}\lg {\cal D},{\cal D}\rg + n d_{\cal E}({\cal D})
- \lg {\cal K}_{\cal E}, {\cal D}\rg + c_1^2+\frac{c_2}{12}
$$
Hence, for $n$ big enough, 
it is positive. By Serre-duality, we know that 
$$
{\rm h}^2({\cal E},{\cal D}')=
{\rm h}^0({\cal E}, {\cal K}_{\cal E}-{\cal D}')=0,
$$
since the latter has a negative intersection with the generic fiber
${\cal E}_b$. Thus, the class of ${\cal D}'$ contains an effective divisor
and ${\cal D}'\cap {\cal E}_b=d_{\cal E}({\cal D})=d_{\cal E}$. 
Then the divisor ${\cal M}$ 
is obtained from this effective divisor by removing
the vertical components (clearly, ${\cal M}$ is irreducible). 

\begin{cor}
The order of $[{\cal E}]\in \H^1(B,{\cal J})$ is equal to 
$d_{\cal E}$.
\end{cor}

\begin{defn}\label{defn:order-m}
A  multisection ${\cal M}$ is said to be torsion {\em of order $m$} 
if $m$ is the smallest positive
integer such that for any $b\in B$ and any pair of points 
$p_b,p_b'\in {\cal M} \cap {\cal E}_b$ the image of the zero-cycle
$p_b-p_b'$  in ${\cal J}_b$ is torsion of order $m$. 
We call a multisection an $nt$-multisection (non-torsion), if 
for a general point $b\in B$ there exist two points 
$p_b,p_b'\in {\cal M} \cap {\cal E}_b$ such that the zero-cycle
$p_b-p_b'\in {\cal J}_b$ is non-torsion.
\end{defn}

\begin{lm}\label{lm:nt}
If an irreducible multisection ${\cal M}\subset {\cal E}$ is not a torsion 
multisection of order $n$ for any $n\in \N$ then 
${\cal M} $ is an $nt$-multisection.
\end{lm}

{\em Proof.} We work over $\C$. 
The union of all torsion multisections 
of ${\cal E}$ is a countable union of divisors. 
So it can't cover all of ${\cal M}$ unless
${\cal M}$ is contained in some torsion multisection.

\

Let ${\cal E}\ra B$ be an elliptic fibration (without multiple fibers). 
There is a natural set of elliptic fibrations  
${\cal J}^m={\cal J}^m({\cal E})$ over
$B$ parametrizing classes of cycles of degree $m$ on the
generic fiber of ${\cal E}$. The
Jacobian of each ${\cal J}^m$ is isomorphic 
to ${\cal J}^0={\cal J}$. The class 
$[{\cal J}^m]\in \H^1_{\rm an}(B,{\cal J})$ equals 
$m\cdot [{\cal E}]$.  
If $\beta$ is a cocycle defining ${\cal J}^{\beta}$ then
the cocycle $m\cdot \beta$ is defined by pointwise multiplication 
by $m$ in the Jacobian fibration. Thus, we obtain
the fibration ${\cal J}^{m\beta}$.

We have natural rational maps 
of algebraic varieties ${\cal J}^m\times_B {\cal J}^k\ra {\cal J}^{m+k}$
which fiberwise is the addition of cycles. These maps provide
the set of isomorphism classes of 
${\cal J}^m$ with the structure of a cyclic group. 
The order of this group for algebraic ${\cal E}$ coincides
with $d_{\cal E}$. 
The identification ${\cal J}^{d_{\cal E}}= {\cal J}^0$ is 
not canonical. It is defined modulo the action of 
$\H^0(B,{\cal J})$ on ${\cal J}^{d_{\cal E}}$ (for example, 
choosing a multisection of degree $d_{\cal E}$ in  
${\cal J}^1={\cal E}$ will fix the identification).

The construction provides  
maps $\eta^m: {\cal J}^1\ra {\cal J}^m$ for any $m\in \N$, since 
${\cal J}^1$ imbeds diagonally into the fiber product 
of $m$ copies of ${\cal J}^1$. 
All the above maps are well defined 
on the open subvarieties
obtained by deleting the singular fibers of the fibrations. 
They are algebraic and they extend to meromorphic maps.

We obtain an action of 
${\cal J}={\cal J}^0$ on ${\cal E}={\cal J}^1$
which is regular in non-singular points of the fibers of 
${\cal J}$ and ${\cal E}$ and which induces a 
transitive action of the fibers ${\cal J}_b$ on ${\cal E}_b$ 
(for smooth fibers).

The maps $\eta^m$ allow to transfer irreducible multisections
between the elliptic fibrations ${\cal J}^m$ (modulo $d_{\cal E}$).
More precisely, we have

\begin{lm}\label{lm:transfer}
Let ${\cal M}\subset {\cal J}^k$ be a torsion multisection 
of order $t$ (with $d_{{\cal J}^k}|t$).
Consider the map $\eta^m\,:\, {\cal J}^k \ra {\cal J}^{km}$.
Then $ \eta^m({\cal M})\subset {\cal J}^{mk}$ 
is a torsion multisection of order exactly $t/{\rm gcd}(t,m)$.
Moreover, if ${\cal M} $ is non-torsion or torsion of order coprime to $m$
then the restriction $\eta^m\,:\, {\cal M}\ra \eta^m({\cal M})$
is a birational map and  
$d_{{\cal J}^k}({\cal M})=d_{{\cal J}^{km}}(\eta^m({\cal M}))$.
\end{lm}

{\em Proof.}
Locally, we have a Jacobian elliptic
fibration and the map $\eta^m$ is multiplication by $m$. 
Therefore, if $\eta^m(x)=\eta^m(y)$ for some $x,y\in {\cal J}^k$
then $m\cdot(x-y)=0$ (in ${\cal J}$). Since ${\cal M}$ is
irreducible, either $(x-y)$ is torsion of order ${\rm gcd}(t,m)$ 
for any pair of points $x,y\in {\cal J}^k_b$ for a general
fiber $b$ or these pairs constitute a divisor in ${\cal M}$.
In the latter case, it follows that the restriction of $\eta^m$ to
${\cal M}$ is a birational map and hence $\eta^m({\cal M})$ is
a multisection of ${\cal J}^{mk}$ of the same degree.

\begin{cor}\label{cor:ep}
Let $p$ be a prime number and ${\cal E}$ an elliptic fibration
with $d_{\cal E}=p$. Let ${\cal M}$ be a torsion multisection 
of ${\cal E}$. Then ${\cal M}$ admits a surjective map onto one
of the $p$-torsion multisections of ${\cal E}$ or onto one
of the non-zero $p$-torsion multisections of ${\cal J}({\cal E})$.
\end{cor}

{\em Proof.}
Suppose that ${\cal M}$ is a torsion multisection of order
$p^kt$ where $(t,p)=1$ and $k\ge 1$. If $k=1$ choose an $\a$ such 
that $\a t = 1 $ mod $p$. 
We have a map
$$
\eta^{\a t }\,:\, {\cal J}^1({\cal E})\ra 
{\cal J}^{\a t}({\cal E})\simeq {\cal J}^1({\cal E}) 
$$
and $\eta^{\a t }({\cal M})$ is a torsion multisection of
order $p$. 

If $k>1$, then $\eta^{p^{k-1}t}({\cal M})$ is a non-trivial 
$p$-torsion multisection (but not a section) in ${\cal J}({\cal E})$.

\begin{defn}\label{defn:tau} 
Let ${\cal Z}$ be any cycle of degree 
$d_{\cal E}({\cal Z})$ on ${\cal E}$
which is given by a combination of 
multisections with integral coefficients. 
We define a class map
$$
\tau_{\cal Z}\, :\,  {\cal E} \ra {\cal J}
$$
by the following rule:
$$
\tau_{\cal Z}(p)=[d_{\cal E}({\cal Z})\cdot p -   
{\rm Tr}_{\cal Z}(\varphi(p))]
$$
for $p\in {\cal E}$. Here we denote by 
${\rm Tr}_{\cal Z}(b)$ the zero-cycle ${\cal Z}\cap {\cal E}_b$.
\end{defn}

\subsection{Monodromy}

Denote by $b_1,...,b_n$ the 
set  of points in $B$ corresponding to singular fibers of ${\cal E}$.
Consider the analytic fibration ${\cal E}^*\ra B^*$, where
$B^*=B\backslash \{b_1,...,b_n\}$, 
obtained by removing all singular fibers from ${\cal E}$. 
We have a natural action of the free group $\pi_1(B^*)$ 
on the integral homology of the fibers. 
The group of automorphisms of the integral homology of 
a generic fiber ${\cal E}_b$ which preserve orientation
is the group $\SL_2(\Z)$. 
One of the main characteristics of
an elliptic fibration ${\cal E}\ra B$ is its
global monodromy group $\Ga$. 

\begin{defn}
The global monodromy group $\Ga=\Ga({\cal E})$ of 
${\cal E}\ra B$ is the image
of $\pi_1(B^*)$ in $\SL_2(\Z)$. 
Denote by $\ind(\Ga)=[\SL_2(\Z):\Ga]$ the index of the
global monodromy. 
A cycle around a point $b_i$ 
(for any $i=1,...,n$) defines a conjugacy class in $\pi_1(B^*)$.
The corresponding conjugacy class in $\SL_2(\Z)$ is called local
monodromy around $b_i$. So we obtain (a class of) cyclic subgroups 
$T_i\subset \SL_2(\Z)$ (up to conjugation). 
\end{defn} 

\begin{rem}
The monodromy group  $\Ga=\Ga({\cal E})$ of an
elliptic fibration ${\cal E}$ coincides with $\Ga({\cal J})$
of the corresponding Jacobian 
elliptic fibration. In particular, for locally isotrivial
elliptic fibrations the monodromy group 
$\Ga({\cal E})$ is a finite subgroup of
$\SL_2(\Z)$. For non-isotrivial 
elliptic fibrations ${\cal E}\ra \P^1$
we have $\ind(\Ga({\cal E}))<\infty$. 
\end{rem}

The group $\SL_2(\Z)$ has a center $\Z/2\Z$ and we shall denote
by $\Ga_c$ the subgroup of 
$\SL_2(\Z)$ obtained by adjoining the center
to $\Ga$. 

\begin{rem}
A generic elliptic fibration ${\cal E}\ra \P^1$ 
has monodromy group $\SL_2(\Z)$. More precisely,
$\SL_2(\Z)$ has two standard nilpotent generators $a,b$. 
Assume that
all singular fibers of ${\cal E}$ are {\em nodal} (rational curves
with one self-intersection). In this case,
we can select a system of vanishing arcs from some points in $B^*$
so that all $T_i$ split into two clusters $I_a$ and $I_b$   
(of equal cardinality) such that $T_i=\lg a\rg$ for $i\in I_a$ 
and $T_i=\lg b\rg $ for $i\in I_b$ (cf. \cite{friedman-morgan-94}, p. 171).
In particular, any two local monodromies corresponding to different
classes generate $\SL_2(\Z)$. 
\end{rem}

Jacobian elliptic fibrations over $\P^1$ arise in families
${\cal F}_r$ parametrized by an integer $r$ 
which is defined through 
the standard Weierstrass form 
$$
y^2=x^3+p(t)x+q(t)
$$
where $p$ (resp. $q$) is a  polynomial of degree $4r$ 
(resp. $6r$), satisfying some genericity conditions. 
There are lists of possible singular fibers,
possible local monodromy groups 
and actions of these groups on the torsion sections of the 
nearby fibers as well as a list
of possible torsion groups of the singular fibers 
(cf. \cite{kodaira} or \cite{barth-peters-vdv-84}).

\

For any non-isotrivial elliptic fibration 
${\cal E}\ra B$ we have a map
$j_B : B\ra \P^1$ defined by $j_B(b):=j({\cal E}_b)$ (where $j$ is the 
standard $j$-invariant of an elliptic curve
with values in $\P^1=\overline{PSL_2(\Z) \backslash H}$). 

\begin{rem}
If $B=\P^1$ and ${\cal E}\in {\cal F}_r$ 
then $j_B=\frac{4p^3}{4p^3+27q^2}$.
Hence, the degree of the map $j_B$ in this case 
is bounded by $12 r$. 
\end{rem}

\begin{prop} \label{prop:mon}
Let ${\cal E}\ra B$ be a non-isotrivial elliptic fibration.
Then  ${\rm ind}(\Ga) \le 2 \deg(j_B)$.
\end{prop}

{\em Proof.} 
The map $j_B$ is the same for 
an elliptic fibration ${\cal E}$ and for the 
Jacobian of ${\cal E}$. Thus we reduce to the
case of Jacobian elliptic fibrations.
Consider the $\Ga$-covering $\tilde{{\cal E}^*}\ra B^*$. It is a
Jacobian elliptic fibration 
over an open analytic curve $\tilde{B^*}$. Since
it is topologically trivial the map $j_B$ lifts to
a holomorphic map 
$\tilde{j}_B : \tilde{B^*}\ra H$  (where $H$ is the upper-half plane). 
This map is $\Ga$-equivariant and it defines
a map $j_{\Ga} : B^*\ra \Ga \backslash H$. Therefore, 
the map $j_B$ on $B^*$ is a composition $j_B= r_{\Ga} \circ j_{\Ga}$,
where $r_{\Ga}$ is the map $r_{\Ga} : 
\Ga \backslash H \ra \SL_2(\Z)\backslash H$. The group $\Ga$ acts
on $H$ through its homomorphism to ${\rm PSL}_2(\Z)$. 
Therefore, the degree of $r_{\Ga}$ is equal to the index 
$\ind(\Gamma_c)$ if $\Gamma$ contains the center 
$\Z/2$ and equal to $\frac{1}{2} \ind(\Gamma)$ otherwise. 

\begin{cor}
The number of possible monodromies in any family of
elliptic fibrations with bounded degree of $j_B$ is 
finite. In particular, for the families ${\cal F}_r$ 
with a given $r$ all monodromy groups have index
$\le 24 r$. 
\end{cor}

\begin{rem}
If we have an algebraic variety which parametrizes
elliptic fibrations then global monodromy changes
only on algebraic subvarieties, where the 
topological type of the projection $\varphi : {\cal E}\ra B$
changes. This variation normally occurs in big codimension.
Indeed, the monodromy is 
completely determined by its action outside of small neighborhoods
of singular fibers. Hence, it doesn't vary 
under small smooth variations of ${\cal E}$. 
\end{rem}

\begin{exam}
For the family ${\cal F}_r$ 
the monodromy is $\SL_2(\Z)$ 
provided that at least two nodal fibers from different
clusters $I_a, I_b$ remain unchanged. The dimension  
of the subvariety in ${\cal F}_r$ with monodromy different from 
$\SL_2(\Z)$ is $\le \frac{1}{2}\dim{\cal F}_r +1$.
\end{exam}

\subsection{Torsion multisections}

In this section we will work over $\C$. 
Let $\varphi :  {\cal E}\ra B$ be an elliptic fibration
and ${\cal M} $ an irreducible multisection of ${\cal E}$. 

\begin{prop}\label{prop:bound-1}
Let ${\cal J}\ra \P^1$ be a non-isotrivial
Jacobian elliptic fibration
with global monodromy $\Gamma\subset \SL(2,\Z)$.
Then there exists a constant $c$ (for example,
$c=\frac{6}{\pi^2}$) such that 
for all torsion multisections
${\cal M}\subset {\cal J}$  of degree $d_{\cal J}({\cal M})$ 
and order $m$ we have
$$
  d_{\cal J}({\cal M})   >\frac{c\cdot m^2}{{\rm ind}({\Gamma})}.
$$
\end{prop}

{\em Proof.} For each $b\in \P^1$ we have an action of 
$\Ga$ on the cycle ${\cal M}\cap {\cal J}_b$ and also an action of
$\Ga $ on the set of points of order $m$ of this fiber.
It follows that ${\cal M}\cap {\cal J}_b$ must coincide with 
an orbit of $\Ga$ on the $m$-torsion points.   
The size of the corresponding orbit for the full
group $\SL(2,\Z)$ acting on {\em primitive} torsion points
of order $m$ (e.g., points of order exactly $m$)
is equal to the product $m^2\cdot\prod_{p|m}(1-1/p^2)$.
Hence, the size of any orbit of $\Gamma$ on the primitive 
$m$-torsion points
of a general fiber is 
$$
>\frac{m^2}{{\rm ind}(\Gamma)}\cdot\prod_{p|m}(1-1/p^2).
$$

\begin{prop}\label{prop:bound}
Let $\Gamma\subset \SL_2(\Z)$ be a subgroup of finite index.
There exists a constant  $m_0(\Ga)$ such that  for all
non-isotrivial Jacobian
elliptic fibrations ${\cal J}\ra B$, with 
at least 4 singular fibers,  
with global monodromy $\Gamma$ and for all 
torsion multisections ${\cal M}\subset {\cal J}$
of order $m$ with $m>m_0(\Ga)$ we have 
$g({\cal M})\ge 2$.
\end{prop}

{\em Proof.} Although this fact is probably well known we 
decided to give an argument. 

Every orbit of the (linear) action of $\Ga$ on $m$-torsion points 
defines an irreducible $m$-torsion 
multisection in ${\cal J}$ (and vice versa). 
Thus we can identify the orbit for a given multisection 
with the quotient  $\Ga/ \Ga'$, 
where $\Ga'$ is a subgroup of finite index in $\Ga$.
The corresponding orbit for a singular fiber ${\cal J}_{b_i}$ 
is equal to the quotient $\Ga/\Ga_i$ where $\Ga_i$ is a subgroup of
$\Ga$ generated by $\Ga'$ and the local subgroup $T_i$ 
(even though $T_i$ are, in principle, 
defined only up to conjugation in
$\Ga$, but specifying the multisection we also specify 
the pair $T_i, \Ga'$ modulo common conjugation).
Therefore, the Euler characteristic of ${\cal M}$ will be
equal to
\begin{equation}\label{eqn:euler}
\chi({\cal M})= |\Ga/\Ga'| \cdot \left( \chi(B) - 
\sum_i (1- \frac{1}{[\Ga'T_i:\Ga']})\right).
\end{equation}
In order to prove our theorem it suffices to observe
that this formula implies 
the growth of the genus of ${\cal M}$ as $m\ra \infty$. 

The fibration ${\cal J}\ra B$ contains at least one
fiber of potentially multiplicative
reduction (pullback of $\infty$ of the $j$-map). 
The local monodromy around any fiber of this type
is an infinite cyclic group which
includes a subgroup of small index (2,3,4,6) generated by
the unipotent transformation 
$\left(\begin{array}{cc}
1 & k \\ 0 & 1 \end{array}\right)$ where $k$ is the number
of components in the fiber.  
The number of $m$-torsion elements in this singular fiber
is at most $m\cdot k$. But the degree of the torsion 
multisection grows like $m^2 $ (cf. \ref{prop:bound}).
Hence $[\Ga'T_i:\Ga']>c'\frac{m}{k\ind(\Ga)}$ 
(for some constant $c'$). 

Similarly, for singular fibers with potentially good reduction
the corresponding local monodromy groups
are among the standard finite subgroups
of $\SL_2(\Z): \Z/2\Z, \Z/3\Z, \Z/4\Z, \Z/6\Z$. 
They all act effectively
on the points of order $m$ for $m\ge 5$. Hence, 
the index $[\Ga'T_i:\Ga']$ (for $m\ge 5$) can be 6,4,3,2. Thus
every fiber of this type contributes at least 
5/6, 3/4, 2/3 or 1/2, respectively. 

Asymptotically, for $m\gg 0$, the contribution from 
every singular fiber of potentially multiplicative 
reduction will tend to 1, the contribution from other
fibers is $\ge 1/2$. Since we have at least 4 singular fibers,
the theorem follows. 

\

\noindent
We have a similar result for non-Jacobian elliptic fibrations:

\begin{prop}\label{prop:genus-estimate}
Let $\Ga\subset \SL_2(\Z)$ be a subgroup of finite index.
There exists a $p_0>0$ (which depends only on $\Ga$) such that
for every  elliptic surface 
${\cal E}\ra\P^1$ with at least 4 singular fibers,
global monodromy $\Ga$ and for
any  torsion multisection ${\cal M} \subset S$
of order $p>p_0$ (where $p$ is a prime number)
we have $g({\cal M})\ge 2$.
\end{prop}

{\em Proof.}
The minimal index of a proper subgroup of
$\SL_2(\Z/p\Z)$ grows with $p$. This implies that
for any subgroup $\Ga$ of finite index in $\SL_2(\Z)$
its projection onto $\SL_2(\Z/p\Z)$ is surjective
for all $p>p_0$. 

Let us prove first that the generic fiber of ${\cal M}$ is
isomorphic to $\Z/p\Z \oplus \Z/p\Z $ (as a $\Ga$-module).
Indeed, the generic fiber corresponds to the orbit
in $\Z/p\Z \oplus \Z/p\Z$ of the action 
of some extension $\tilde{\Ga}$ of $\SL_2(\Z/p\Z)$ by a subgroup 
of translations.
If the latter contains a non-trivial translation, then 
it contains the whole module of translations $\Z/p\Z \oplus \Z/p\Z$,
hence the claim. Otherwise, $\tilde{\Ga}$ is isomorphic to 
$\SL_2(\Z/p\Z)$ and it acts on $\Z/p\Z \oplus \Z/p\Z$ linearly.
Indeed, it contains a central involution which has
exactly one invariant element in $\Z/p\Z \oplus \Z/p\Z$.
This linearizes the action which leads to a contradition 
if  ${\cal E} $ is non-Jacobian. If ${\cal E} $ is a Jacobian
elliptic fibration then the degree $d_{\cal E}({\cal M})=p^2-1$.

Now we can estimate from below the genus of ${\cal M}$ using 
a formula similar to \ref{eqn:euler}.
More precisely, we have
$$
\chi({\cal M})\le -\frac{p^2-1}{2} -c,
$$
where $c$ is a (small, effectively computable) constant.
 
\begin{prop}
\label{prop:genus-estimate-2}
For every subgroup $\Ga \subset \SL_2(\Z)$ of finite index  
there exists a $p_0$ such that for all primes $p>p_0$ and all
(non-isotrivial) non-Jacobian  elliptic fibrations ${\cal E}\ra\P^1$
of degree $d_{\cal E}=p$,  
with at least 4 singular fibers
and with global  monodromy $\Ga$ 
every torsion multisection ${\cal M}$ of ${\cal E}$ 
has genus $g({\cal M})\ge 2$. 
\end{prop}

{\em Proof.}
First observe that the class of order $p$ in the Shafarevich
group corresponds to a cocycle with coefficients in the
$p$-torsion sub-sheaf of ${\cal J}({\cal E})$. 
Therefore, the elliptic fibration ${\cal E}$ corresponding 
to a cocycle of order $p$ contains a $p$-torsion 
multisection ${\cal M}$. 

By \ref{cor:ep}, we know that every torsion multisection 
${\cal M}'\subset {\cal E}$ 
admits a map onto the $p$-torsion multisection in 
${\cal M}\subset{\cal E}$
or a $p$-torsion multisection in the corresponding Jacobian 
elliptic fibration ${\cal J}({\cal E})$. 
Now we apply \ref{prop:genus-estimate}.

\begin{prop}\label{prop:no-rat}
Let ${\cal E}\ra \P^1$ be an elliptic fibration 
(with at least 4 singular fibers and 
fixed monodromy group $\Ga$ as above).
Let $p>p_0$ a prime
number not dividing  the degree $d_{\cal E}$. 
Let ${\cal E}'\ra \P^1$ be an elliptic fibration 
of degree $p\cdot d_{\cal E}$, obtained by
dividing the cocycle corresponding to ${\cal E}$ by $p$.
Then ${\cal E}'$ has no rational or 
elliptic torsion multisections. 
\end{prop}

{\em Proof.}
Let ${\cal E}''={{\cal E}'}^{d_{\cal E}}$.
It is a fibration of order $p$ (with the same monodromy
group $\Ga$.)
Any torsion multisection of ${\cal E}'$ is mapped to 
a torsion multisection of ${\cal E}''$.
By \ref{prop:genus-estimate-2}, the genus of
any torsion multisection in ${\cal E}''$, and therefore in ${\cal E}'$
is $\ge 2$.

\begin{lm}\label{lem:4-fibers}
Any elliptic K3 surface 
$S\ra \P^1$ with $\Pic(S)\le 19$ has at least
4 singular fibers, including at least one potentially
multiplicative fiber. 
\end{lm}

{\em Proof.} 
The proof is topological and works for Jacobian and 
non-Jacobian elliptic fibrations.
Denote by $\chi({\cal E}_b)$ the Euler characteristic 
and by $r({\cal E}_b)$ the rank of the lattice spanned by classes of the
irreducible components of the singular fiber ${\cal E}_b$. 
Then $\chi({\cal E}_b) -r({\cal E}_b)=1$ if the fiber has multiplicative 
reduction (Type $I_n$), or $\chi({\cal E}_b)-r({\cal E}_b)=2$ otherwise. 
We have $\sum \chi({\cal E}_b) =24$ and 
$\sum  r({\cal E}_b)\le 18$ (for more details see, 
for example \cite{ye}, pp. 7--9).

\begin{rem}
In \cite{beauville-82-2} 
Beauville proves that every {\em semi-stable}
non-isotrivial elliptic fibration has at least 4 singular fibers
and classifies those which have exactly 4. (These are 6 
modular families, cf.  \cite{beauville-2}, p. 658.)  
There is a complete classification of elliptic K3 surfaces with 
3 singular fibers in \cite{ye}.
For recent work concerning the minimal
number of singular fibers in fibrations with generic fiber
a curve of genus $\ge 1$ see \cite{tan}, \cite{ye}. 
\end{rem}

As a corollary we obtain Theorem 
\ref{thm:nt-inf} stated in the introduction:

\begin{cor}\label{coro:k3-rat}
Every algebraic elliptic K3 surface 
$S\ra \P^1$ with $\rk \, \Pic(S)\le 19$
has infinitely many
rational $nt$-multisections.
\end{cor}

{\em Proof.}
If $S$ is Jacobian we denote
by $S'$ some algebraic non-Jacobian 
elliptic K3 surface with Jacobian ${\cal J}(S')=S$.
Otherwise, we put $S'=S$.  
Dividing (the cocycle defining) $S'$ by different primes $p>p_0$
we obtain elliptic K3 surfaces ${\cal E}_p$ (of different degrees). 
By proposition \ref{prop:no-rat}, ${\cal E}_p$ don't 
contain rational or elliptic torsion multisections. At the same
time, by deformation theory, they contain rational multisections
of degree divisible by $d_{{\cal E}_p}$. Therefore, 
we can produce a sequence of rational $nt$-multisections in $S'$
(and consequently, in $S$) 
of increasing degrees.

\section{Density of rational points}

\subsection{Multisections}

From now on we will work over a number field $K$ and we  
restict to the case of the base $B=\P^1$.

\begin{prop}\cite{bogomolov-tschi-98}\label{theo:non-torsion}
Let $ \varphi_{\cal E}  :  {\cal E}\ra \P^1$ be an elliptic
fibration defined over $K$ with a $nt$-multisection ${\cal M}$.
Then for all but finitely many 
$b\in \varphi_{\cal J}({\cal M}(K))\subset \P^1(K)$ 
the fibers ${\cal E}_{b}$ have infinitely many rational points. 
\end{prop}

{\em Proof.} Since ${\cal M}$ is an $nt$-multisection, we have a 
birational map 
$$
\tau\,:\, {\cal M}\ra \tau ({\cal M})\subset {\cal J}({\cal E}).
$$
An argument using
Merel's theorem (or simply base change to $\tau({\cal M})$) implies
that rational points are dense in the fibers ${\cal J}_b$ for almost all 
$b\in \varphi_{\cal J}(\tau({\cal M})(K))$
(for a sufficiently large finite extension of $K/\Q$). 
Then one can translate points in $({\cal E}_b\cap {\cal M})(K)$ 
(for $b\in \varphi({\cal M}(K))$) to 
obtain a Zariski dense set of rational points in the fibers ${\cal E}_b$ and 
consequently in ${\cal E}$.  

\begin{cor} Let $S\ra \P^1$ be an elliptic K3 surface
defined over a number field $K$. Then rational points
on $S$ are potentially dense.
\end{cor}

{\em Proof.}
By \ref{coro:k3-rat}, every algebraic elliptic K3 surface 
with $\rk \, \Pic(S)\le 19$ has infinitely many
rational $nt$-multisections. If $\rk \, \Pic(S)=20$ 
we use \ref{prop:auto}.

\begin{defn}\label{dfn:non-torsion-section}
Let $\varphi  : {\cal E}\ra B$ 
be an elliptic fibration. 
A {\em saliently ramified multisection} of ${\cal E}$ is 
a multisection ${\cal M}$ 
which intersects a fiber ${\cal E}_b$
at some smooth point $p_b$ with local intersection multiplicity $\ge 2$.
\end{defn}

\begin{prop}\label{prop:salient-non-torsion}\cite{bogomolov-tschi-98-2}
Suppose that ${\cal M}\subset {\cal E}$ is a saliently ramified rational or
elliptic multisection. Then it is an $nt$-multisection. 
Consequently, rational points on ${\cal E}$ are potentially dense. 
\end{prop}

\begin{cor}
Let $S$ be an algebraic surface admitting two elliptic fibrations over
$\P^1$.  Then rational points on $S$ are potentially dense.
\end{cor}

\begin{rem}
An alternative approach to potential density of rational points on 
elliptic K3 surfaces 
${\cal E}\ra \P^1$  would be to show that there exists an family of
elliptic curves ``transversal'' to the given elliptic fibration. 
Then a generic elliptic curve in the transversal elliptic fibration is
a saliently ramified multisection of ${\cal E}\ra \P^1$. It remains to apply
\ref{prop:salient-non-torsion}. 
\end{rem}

\subsection{Automorphisms}

Let $X$ be a K3 surface defined over a number field $K$.
We have a hyperbolic lattice $\Pic(X):= \Pic(X_{\C})\subset {\cal L}$ 
where ${\cal L}= 3\cdot {\cal H}\oplus 2\cdot (-E_8)$ and a monoid
of effective divisors $\L(X)\subset \Pic(X)$. We denote
by $\Aut(X)$ the group of (regular) 
algebraic automorphisms of $X$ (over $\C$). 
Observe that $\Aut(X)$ is finitely generated. 
We can guarantee that $\Aut(X)$ is defined over $K'$, for some finite 
extension $K'/K$.  

\begin{rem} V. Nikulin proved that there are only finitely 
many isomorphy types of lattices $\Pic(X)$ for K3 surfaces
with $\rk\,\Pic(X)\ge 3 $ such that the corresponding group
$\Aut(X)$ is finite (cf.  \cite{nikulin-86}).  
We can prove potential density for those surfaces 
from Nikulin's list which contain (semipositive)
elements with square zero.
For example, there are  17 lattices that give finite 
automorphism groups $\Aut(X)$ for $\rk\, \Pic(X)=4$ and of those
17 lattices 15 contain elements with square zero (and therefore 
admit elliptic fibrations) (cf. \cite{vinberg},
\cite{nikulin-86}). 
\end{rem}

\begin{exam}
There exists a K3 surface of rank 4 with 
the following Picard lattice:
$$
\left(
\begin{array}{rrrr}
2  & -1 & -1 & -1  \\
-1 & -2 &  0 & 0   \\
-1 & 0  & -2 & 0   \\
-1 & 0  &  0 &-2   \end{array}
\right)
$$
There are no elements or square zero and the group of automorphisms
$\Aut(X)$ is finite. We don't know whether or not rational points on 
$X$ are potentially dense.  
\end{exam}

\begin{lm} Suppose that $\Aut(X)$ is infinite.
Then $\L(X)$ is not finitely generated.
\end{lm}

{\em Proof.} If suffices to identify $\Aut(X)$ (up to a finite index)
with the subgroup of $\Aut({\cal L})$ which preserves $\L(X)$. 
The set of generators of $\L(X)$ is preserved under $\Aut(X)$. 
If this set is finite $\Aut(X)$ must be finite as well. 

\begin{thm} \label{prop:auto}
Let $X$ be a K3 surface over a number field $K$ 
with an infinite group of automorphisms. 
Then rational points on $X$ are potentially dense. 
\end{thm}

{\em Proof.}
It suffices to find a rational curve $C\subset X$ 
such that the orbit of $C$ under $\Aut(X)$ is infinite.
The monoid $\L(X)$ is generated by classes of $(-2)$-curves,
curves with square zero and primitive classes with positive
square. It follows from (\ref{rem:prim}) that every
generator of $\L(X)$ is represented
by  a (possibly singular) irreducible rational curve.  
Suppose that orbits of $\Aut(X)$ on the generators of $\L(X)$ 
are all finite. Then the group $\Aut(X)$ is finite and the number
of elements is bounded by a function depending only 
on the rank of the lattice. (Indeed, any group acting
on a lattice of rank $n$ embedds into $\SL_n(\Z_3)$. The 
normal subgroup of elements in $\SL_n(\Z_3)$  equal to the identity
modulo 3 consists of elements of infinite order. Hence any
subgroup of the automorphisms of the lattice has a subgroup 
of finite index which consists of elements of infinite order.)
So there exists an element of infinite order. For this element
the orbit of some generator of $\L(X)$ is infinite. This class
is represented by a rational curve $C$. The orbit of $C$ is
not contained in any divisor in $X$. 
Extending the field, if 
necessary, we can assume that rational points on $C$ are Zariski
dense. This concludes the proof.

\begin{rem}
Certainly, there are algebraic varieties $X$ such that the orbit 
under $\Aut(X)$ of any  given rational point is 
always contained in a divisor. For example, consider a generic
Jacobian elliptic surface ${\cal J}$
with a non-torsion group of sections.
Then $\Aut({\cal J})$ is generated by the group 
of fiberwise involutions with respect to the sections. 
In particular, inspite of the fact that 
the group is infinite the fibers are preserved and the orbit of 
any point is contained in a divisor. 
However, rational points on $X$ are Zariski dense, as there is
a rational section of infinite order (in $\Aut(X)$ and in ${\cal J}$). 
\end{rem}

\begin{cor} \label{cor:}
Let $X$ be a K3 surface such that $\rk \, \Pic(X)\ge 2$ and 
$\Pic(X)$ contains no classes with square zero and square $(-2)$.
Then $\Aut(X)$ is infinite and rational points on $X$ are
potentially dense. 
\end{cor}

{\em Proof.} The monoid $\L(X)$ is infinitely generated.

\bibliographystyle{plain}
\bibliography{acl}

\end{document}